\documentclass[runningheads]{llncs}
\usepackage[T1]{fontenc}
\usepackage{newtxtext}       %
\usepackage[varvw]{newtxmath}       

\usepackage{graphicx}
\usepackage{amsmath,mathtools}
\usepackage{bm,mathrsfs,cancel,stmaryrd}
\usepackage{systeme}
\usepackage{booktabs,float,enumitem}
\usepackage{xcolor}
\usepackage{siunitx}
\usepackage{import}
\usepackage[caption=false,font=small]{subfig}

\usepackage{tikz}
\usetikzlibrary{arrows.meta,positioning}
\tikzset{
    ddnode/.style={circle,draw,thick,fill=white,inner sep=2pt,minimum size=7pt},
    arcpos/.style={->,thick,>=Latex},
    arcneg/.style={arcpos,dashed},
    lab/.style={font=\small}
}

\newcommand\R{\ensuremath{\mathbb{R}}}

\DeclareMathOperator{\conv}{conv}
\DeclareMathOperator{\proj}{Proj}

\DeclareMathOperator{\rank}{rank}

\let\svlim\lim\def\lim{\svlim\limits}

\let\epsilon\varepsilon

\definecolor{correct}{HTML}{009900}

\newcommand\hide[1]{}

\spnewtheorem{assumption}{Assumption}{\bfseries}{\itshape}

\makeatletter
\newcommand{\pushright}[1]{\ifmeasuring@#1\else\omit\hfill$\displaystyle#1$\fi\ignorespaces}
\newcommand{\pushleft}[1]{\ifmeasuring@#1\else\omit$\displaystyle#1$\hfill\fi\ignorespaces}
\makeatother

\usepackage[hidelinks]{hyperref}

\newcommand{\X}{\mathcal{X}}
\newcommand{\Y}{\mathcal{Y}}

\newcommand{\F}{\Omega}
\newcommand{\setR}{\mathcal{R}}

\newenvironment{proofof}[1]
  {\renewcommand{\proofname}{Proof of #1}\begin{proof}}
  {\end{proof}}
  
\begin{document}
\title{Valid Inequalities for Bilevel 0-1 Shortest-Path Problems}
\titlerunning{Valid Inequalities for Bilevel 0-1 Shortest-Path Problems}
%
\author{Bruno M. Pacheco, Margarida Carvalho, and Kim Yu}
\authorrunning{B. M. Pacheco, M. Carvalho, and K. Yu}
%
\institute{CIRRELT, Université de Montréal}
%
\maketitle              
\begin{abstract}
    We characterize valid inequalities for bilevel programs in which the leader makes binary decisions, and the follower solves a shortest-path problem.
    We exploit a dualize-and-combine reformulation and show that even with ideal formulations for the leader and the follower primal solution sets, fractional solutions may still appear in the continuous relaxation.
    From that reformulation, we identify that the interaction between the leader variables and the follower dual variables, which we capture through a leader--potential set, is the key structure to be convexified to guarantee that a continuous relaxation recovers the optimal value of the original problem, closing the objective gap.
    Furthermore, we characterize the valid inequalities for the leader--potential set through network flows, and identify a family of inequalities we name \emph{root cuts} that suffice to recover the optimal value. Complementary to the structural results, we prove that the cut generation problem is NP-hard even with partially fixed coefficients and a network with two nodes and two parallel arcs.
    Motivated by the complexity result, we propose a constructive method to fix coefficients in the valid inequalities while preserving separation guarantees, by exploiting the network structure.

\end{abstract}

\section{Introduction}

Bilevel programs model hierarchical decisions in which a \emph{leader}, at the \emph{upper level}, anticipates the optimal response of a \emph{follower}, at the \emph{lower level}~\cite{dempeBilevelOptimizationTheory2020,kleinertSurveyMixedIntegerProgramming2021}.
In this work, we are interested in bilevel binary programs of the form
\begin{equation}\label{eq:bilevel-problem}
\begin{aligned}
    \min_{\bm{x},\bm{y}}\quad & \bm{c}^{\top}(\bm{x},\bm{y})\\
    \text{s.t.}\quad & \bm{x}\in \X,\\
    & \bm{y}\in\arg\min_{\bm{y}'\in \Y}\bm{q}(\bm{x})^{\top}\bm{y}'.
\end{aligned}
\end{equation}
Here, $\bm{x}\in \X\subseteq\{0,1\}^{n}$ represents the leader's decisions, where $\X$ is a set of binary vectors described by linear inequalities. Vector $\bm{y}$ represents the follower's decisions, and $\bm{c}$ collects the coefficients of the leader's linear objective in $(\bm{x},\bm{y})$.
The follower solves a shortest-path problem, and $\Y$ consists of the arc-incidence vectors of the feasible paths.
The leader affects the follower's optimization solely through arc costs $\bm{q}(\bm{x})$ that are affine in $\bm{x}$.
The formulation adopts the optimistic assumption, under which ties among optimal follower paths are resolved in the leader's favor.

Shortest-path followers as above arise in a wide range of optimization settings, and in many of these contexts, the leader naturally affects the follower's objective.
Network pricing is one such example, where a leader anticipates how users (followers) will route through a network in response to her decisions~\cite{labbeBilevelModelTaxation1998,brotcorne2001bilevel,buiCatalogFormulationsNetwork2022}. 
The same structure also arises when a combinatorial follower is represented by a finite-stage dynamic program~\cite{buiSolvingCombinatorialPricing2025}, such as a decision diagram~\cite{lozanoConstrainedShortestPathReformulations2024,bergmanDecisionDiagramsOptimization2015}.
In these representations, which may be exponentially large, paths through an acyclic state-transition graph encode follower solutions, and arc weights accumulate their objective values.
When the follower's feasible set is fixed, and its objective coefficients depend affinely on binary leader decisions, this yields precisely the follower structure in~\eqref{eq:bilevel-problem}.
Similarly, our setting accommodates certain interdiction models after transferring linking constraints---lower-level constraints parametric on the leader's variables---to the follower's objective and constructing a network reformulation based on the follower's dynamic-programming representation. Knapsack interdiction provides such an example~\cite{capraraBilevelKnapsackInterdiction2016,fischettiInterdiction2019,buiSolvingCombinatorialPricing2025}.
Related network formulations encode leader-dependent follower feasibility through arc capacities or represent the follower's value function over leader decisions~\cite{vasquezSinglelevelReformulationBinary2025,lozanoNetworkRelaxationsCombinatorial2026}.

Under common assumptions, the standard arc-flow formulation of the follower's shortest path is ideal. Thus, traditional \emph{dualize-and-combine} can be directly applied to solve the bilevel problem, enforcing follower's optimality through primal feasibility, dual feasibility, and strong duality~\cite{buiCatalogFormulationsNetwork2022,lozanoConstrainedShortestPathReformulations2024,buiSolvingCombinatorialPricing2025}.
In our setting, this yields an extended single-level formulation in the leader decisions, network paths (follower's primal), and node potentials (follower's dual),
linked by a bilinear optimality constraint.
However, even with ideal descriptions of the leader and follower feasible sets, the resulting single-level relaxation can contain fractional points outside the convex hull of bilevel-feasible solutions.
In other words, the relaxation may still contain fractional optimal solutions that satisfy follower's optimality and have strictly better objective values than any feasible solution to the original problem.

While cutting-planes have demonstrated success in separating fractional solutions in general mixed-integer programming, in bilevel settings we observe an absence of methods targeting fractional points that satisfy the follower-optimality condition.
Existing cut generation for mixed-integer bilevel optimization targets, in particular, points that satisfy all leader and follower constraints, including integrality, but violate follower optimality~\cite{fischettiGeneralPurpose2017,fischettiIntersectionCuts2018}.
When a point with a fractional leader decision already satisfies the follower-optimality condition, Tahernejad and Ralphs~\cite{tahernejadValidInequalitiesMixed2025} discuss the use of standard mixed-integer linear programming cuts.
Recent work by V\'asquez et al.~\cite{vasquezNetworkFlowBasedCuts2026} derives cuts in the original variable space from networks that jointly encode the binary leader assignments and the follower's value function, and strengthens them via network-based tilting.
In contrast, our bilevel problem~\eqref{eq:bilevel-problem} contains a network structure only in the follower, and we explore how this structure can be leveraged to characterize the convex hull of feasible solutions.
Davarnia et al.~\cite{khademniaConvexificationBilinearTerms2025} convexify bilinear terms over a dual network polytope and a simplex leader to strengthen dual formulations of network interdiction.
In our setting, the leader decisions interact with the follower's dual through the arc costs, and not through bilinear terms.
Furthermore, we do not assume the leader to be restricted to simplex variables.

\paragraph{Contributions and organization.}
We describe our problem setting in Section \ref{sec:problem-definition}, and we make the following contributions in Sections \ref{sec:valid-inequalities}--\ref{sec:cut-generation}. 
\begin{enumerate}[label=(\roman*)]
    \item We establish the role of leader--potential convexification in the full single-level formulation (Section~\ref{sec:dual-network}).
    With the follower's path polytope and the bilinear optimality constraint retained, we show that convexifying the leader--potential set makes every remaining feasible point a convex combination of feasible single-level solutions.

    \item We characterize all valid linear inequalities to the leader--potential set (Section~\ref{sec:inequalities-flows}). Specifically, we provide necessary and sufficient conditions for the coefficients of node potential variables through network flow.

    \item For a subclass of valid inequalities with special structures, which we call \emph{root cuts}, we prove that adding them to the continuous relaxation already suffices to recover the optimal value of the original problem (Section~\ref{sec:root-supplies}).

    \item We propose a constructive coefficient fixing for the node potentials enhanced by the network structure (Section~\ref{sec:supply-selection}), and we prove that the cut generation problem is NP-hard even with fixed node potential coefficients and a network consisting of two parallel arcs (Section~\ref{sec:separation-complexity}). 
\end{enumerate}

\section{Problem Setting}\label{sec:problem-definition}

We first detail the follower's problem in Section~\ref{subsec:network}
and then derive the single-level reformulation
of~\eqref{eq:bilevel-problem} via dualize-and-combine
in Section~\ref{subsec:dualizeandcombine}.

\subsection{Network assumptions and path formulation}\label{subsec:network}

Let $\mathcal{D}=(\mathcal{V},\mathcal{A})$ be the directed graph for the follower's problem, with distinct root $r\in\mathcal{V}$ and destination $t\in\mathcal{V}$.
We also assume $\mathcal{D}$ to be acyclic, we allow parallel arcs and retain only nodes and arcs on $r\leadsto t$ paths.
This preprocessing preserves all follower choices and makes $r$ and $t$ the unique source and sink, respectively.
Two network examples are illustrated in Figure~\ref{fig:example-networks}.

The above graph assumptions are not stringent; they naturally accommodate the representations from finite-stage dynamic programs and decision diagrams~\cite{buiSolvingCombinatorialPricing2025,bergmanDecisionDiagramsOptimization2015}, but also broader routing settings through equivalent representations.
When we have multiple independent shortest-path followers, connecting copies of their networks in series yields a single path problem of size linear in their original size, with each path selecting one route per follower.
For cyclic networks, acyclic representations preserve all optimal follower paths if every directed cycle has strictly positive cost for every $\bm{x}\in \mathcal X$.
Both constructions preserve the leader's linear objective and optimistic tie-breaking under projection.
Common assumptions such as positive base costs and nonnegative tolls~\cite{buiCatalogFormulationsNetwork2022} or absence of nonpositive cost cycles~\cite{labbeBilevelModelTaxation1998,brotcorneJointDesignPricing2008} satisfy the condition for an acyclic representation.

The follower's path set is independent of $\bm{x}$, and its arc weights are given by $\bm{q}(\bm{x})=Q\bm{x}+\bm{q}_0$, where $Q\in\R^{\mathcal{A}\times n}$ and $\bm{q}_0\in\R^{\mathcal{A}}$.
For each arc $a\in\mathcal{A}$, let $u(a)$ and $v(a)$ denote its tail and head.
The node--arc incidence matrix $D\in\R^{\mathcal{V}\times\mathcal{A}}$ has column $D_a=\bm{e}_{u(a)}-\bm{e}_{v(a)}$, where $\bm{e}_u$ is the unit vector indexed by node $u$.

Recall that $\mathcal{Y}$ is the set of arc-incidence vectors of $r\leadsto t$ paths.
By flow decomposition and acyclicity, every nonnegative unit flow is a convex combination of path-incidence vectors.
Thus, the follower has a bounded ideal linear formulation and a finite optimal value for every leader assignment.
In particular,
\begin{equation}\label{eq:path-polytope}
    \conv(\mathcal Y)=\{\bm{y}\in\R_+^{\mathcal{A}}:D\bm{y}=\bm{e}_r-\bm{e}_t\}.
\end{equation}

\subsection{Strong-duality reformulation}\label{subsec:dualizeandcombine}

In what follows, let 
\[\F := \left\{(\bm{x}, \bm{\theta}, \bm{y})\in \mathcal X\times\R^{\mathcal{V}} \times \mathcal Y: \; \bm{q}(\bm{x})^{\top}\bm{y}\le\theta_r, \quad D^{\top}\bm{\theta}\le\bm{q}(\bm{x}),\quad \theta_t=0\right\}.\]
The constraints $\bm{x}\in\mathcal X$ and $\bm{y}\in\mathcal Y$ describe leader feasibility and follower path-incidence vectors, respectively.
The \emph{dual arc inequalities} $D^{\top}\bm{\theta}\le\bm{q}(\bm{x})$ couple leader decisions and node potentials.
The bilinear \emph{optimality constraint} $\bm{q}(\bm{x})^{\top}\bm{y}\le\theta_r$ enforces strong duality and couples all three variable blocks, $\bm{x}$, $\bm{\theta}$ and $\bm{y}$. Here, we normalize the dual node potentials $\bm{\theta}$ by $\theta_t=0$. 
Following from \eqref{eq:path-polytope} and strong duality, \eqref{eq:bilevel-problem} admits the following equivalent single-level formulation:

\begin{equation}\label{eq:single-level-reformulation}
\begin{aligned}
    \min_{(\bm{x},\bm{\theta},\bm{y})\in\F}\quad
        & \bm{c}^{\top}(\bm{x},\bm{y}). 
\end{aligned}
\end{equation}
This formulation is exact for binary $(\bm{x},\bm{y})$~\cite{lozanoConstrainedShortestPathReformulations2024}, but its continuous relaxation may have fractional extreme points.
We seek valid inequalities to strengthen it.

\section{Valid Inequalities}\label{sec:valid-inequalities}

This section identifies which valid inequalities must be added to the relaxation of~\eqref{eq:single-level-reformulation} to recover the bilevel optimum.
Throughout this section, we assume $\X\ne\varnothing$, since otherwise \eqref{eq:bilevel-problem} is infeasible and every inequality is valid.
Section~\ref{sec:dual-network} shows that it suffices to convexify the set coupling leader decisions and node potentials.
Section~\ref{sec:inequalities-flows} characterizes its valid inequalities through network supplies and affine upper bounds of minimum-cost flow values over $\mathcal X$.
Section~\ref{sec:root-supplies} shows that the root supply $\bm{e}_r-\bm{e}_t$ suffices for the bilevel objective and establishes the limits of more general unit supplies for full convexification.

We first justify seeking a linear description of the feasible hull despite the bilinear optimality constraint.
Although the extended formulation \eqref{eq:single-level-reformulation} introduces continuous variables, the convex hull of its feasible solutions $\conv\!\left( \F \right)$ is a polytope.
This can be seen from a disjunctive programming perspective, which we formalize in the following proposition.

\begin{proposition}\label{prop:SP-is-polyhedral}
    $\conv(\F)$ is a polytope.
\end{proposition}

Consequently, finitely many valid linear inequalities suffice to recover the bilevel optimum from the continuous relaxation of \eqref{eq:single-level-reformulation}.
To identify a structured source of such inequalities, we first account for leader and follower integrality separately.
The follower's flow formulation already describes $\conv(\mathcal Y)$, and we may also replace $\mathcal X$ by $\conv(\mathcal X)$.
As the next example shows, however, the interaction between the variable blocks can still leave fractional points outside $\conv(\F)$ and an objective gap.

\begin{figure}[htbp]
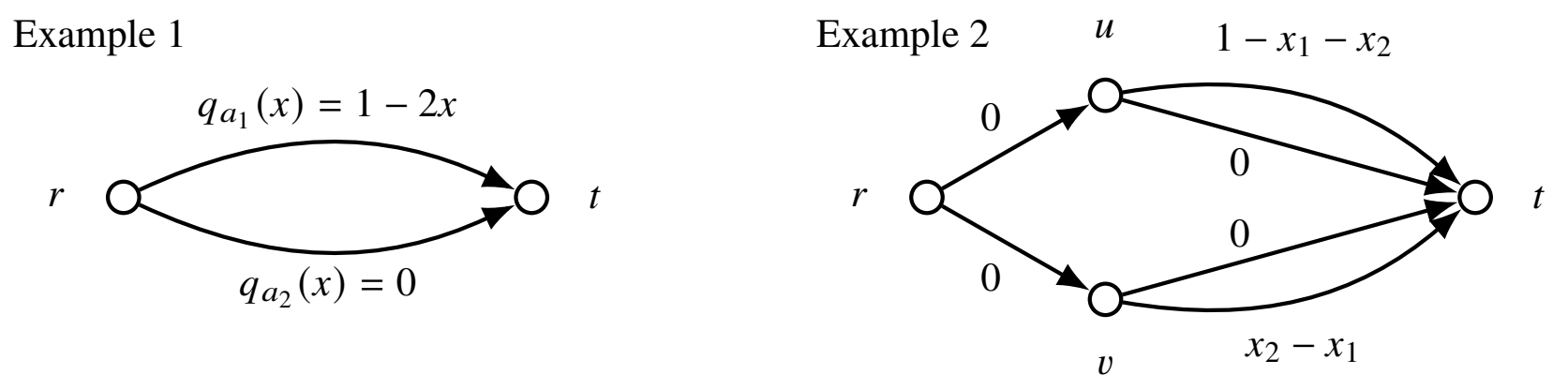

    \centering
            \centering
            \begin{tikzpicture}[scale=.9,transform shape,font=\small]
                \begin{scope}
                    \input{figures/example-1-network.pic.tex}
                    \node[anchor=east,inner sep=0pt] at (.5,1.25)
                        {Example~\ref{eg:full-formulation-fractionality}};
                \end{scope}
            \begin{scope}[xshift=6.3cm]
                \input{figures/example-2-network.pic.tex}
                \node[anchor=east,inner sep=0pt] at (.5,1.25)
                    {Example~\ref{eg:node-pair-cuts-insufficient}};
            \end{scope}
            \end{tikzpicture}
    \caption{Follower networks. Arc labels give costs.}
    \label{fig:example-networks}
\end{figure}

\begin{example}\label{eg:full-formulation-fractionality}
    Let $\mathcal X=\{0,1\}$ and consider two parallel arcs $a_1,a_2$ from $r$ to $t$, with weights $q_{a_1}(x)=1-2x$ and $q_{a_2}(x)=0$.
    The network is illustrated in Figure~\ref{fig:example-networks}.
    The follower's path set is $\mathcal Y=\{(1,0),(0,1)\}$, with $y_i$ indicating arc $a_i$.
    Since $\theta_t=0$ and $y_2=1-y_1$, we retain only $(x,\theta_r,y_1)$.
    Replacing $\mathcal X$ and $\mathcal Y$ by their convex hulls in \eqref{eq:single-level-reformulation} gives $(1-2x)y_1\le\theta_r\le\min\{1-2x,0\}$ for $0\le x,y_1\le1$.
    The fractional point $(1/2,0,1)$ is feasible, but lies outside the segment joining the only binary-feasible solutions $(0,0,0)$ and $(1,-1,1)$.
    For the objective $x-y_1$, it has value $-1/2$, whereas both binary solutions have value zero.
\end{example}

Let $\phi(\bm{x})$ denote the follower's optimal value for $\bm{x}\in\conv(\mathcal X)$.
In the example, the dual arc inequalities impose $\theta_r\le\phi(x)=\min\{1-2x,0\}$, while convexifying the leader--potential set gives the stronger inequality $\theta_r\le-x$.
This motivates studying that set and the effect of its convexification on the full formulation.

\subsection{The dual network polyhedron}\label{sec:dual-network}

For fixed leader decisions, the dual arc inequalities define a dual network flow polyhedron~\cite{borgwardtQuadraticDiameterBounds2016}.
Allowing the arc weights to depend on the leader gives the mixed-integer set
\begin{equation*}\label{eq:S}
    \mathcal S=\left\{(\bm{x},\bm{\theta})\in \mathcal X\times\R^{\mathcal{V}}:
           D^{\top}\bm{\theta}\le\bm{q}(\bm{x}),\ \theta_t=0\right\}
\end{equation*}
and its relaxation
\begin{equation*}\label{eq:P}
    \mathcal P=\left\{(\bm{x},\bm{\theta})\in\conv(\mathcal X)\times\R^{\mathcal{V}}:
           D^{\top}\bm{\theta}\le\bm{q}(\bm{x}),\ \theta_t=0\right\}.
\end{equation*}
We refer to $\mathcal S$ as the \emph{leader--potential} set, and to $\mathcal P$ as its continuous relaxation.
We note that structures similar to $\mathcal{P}$ have been explored under the names \emph{dual network polytope}~\cite{khademniaConvexificationBilinearTerms2025} and \emph{dual network flow polyhedra}~\cite{borgwardtQuadraticDiameterBounds2016} in different settings and with objectives distinct from ours.

Using $\conv(\mathcal X)$ in $\mathcal P$ isolates the interaction between leader decisions and potentials.
For a set $\mathcal C$ in $(\bm{x},\bm{\theta})$-space, write $\mathcal C(\bar{\bm{x}})=\{\bm{\theta}:(\bar{\bm{x}},\bm{\theta})\in\mathcal C\}$ for its \emph{potential fiber}.
Every fiber $\mathcal P(\bar{\bm{x}})$ with $\bar{\bm{x}}\in\conv(\mathcal X)$ is nonempty, as witnessed by shortest-path costs from each node to $t$.
Moreover, $\mathcal S(\hat{\bm{x}})=\mathcal P(\hat{\bm{x}})$ for $\hat{\bm{x}}\in\mathcal X$, so $\proj_{\bm{x}}\mathcal S=\mathcal X$.
The finitely many nonempty fibers of $\mathcal S$ are pointed polyhedra with a common recession cone, hence $\conv(\mathcal S)$ is also a polyhedron.

With the bilinear optimality constraint retained, convexifying $\mathcal S$ excludes every point outside $\conv(\F)$.
Define
\begin{equation*}\label{eq:exact-hull-relaxation}
    \mathcal E=\left\{(\bm{x},\bm{\theta},\bm{y}):
       (\bm{x},\bm{\theta})\in\conv(\mathcal S),\ \bm{y}\in\conv(\mathcal Y),\
       \bm{q}(\bm{x})^{\top}\bm{y}\le\theta_r\right\}.
\end{equation*}

\begin{proposition}\label{thm:exact-convexification}
    $\conv(\mathcal E)=\conv(\F)$.
\end{proposition}

Thus every point of $\mathcal E$ is a convex combination of feasible single-level solutions, and every linear objective has the same optimal value over $\mathcal E$ and $\F$.
This justifies focusing on valid inequalities for $\mathcal S$.
Although $\mathcal E$ may be nonconvex, these inequalities also strengthen the McCormick relaxation of the bilinear terms.
In \autoref{eg:full-formulation-fractionality}, adding $\conv(\mathcal S)=\{(x,\theta_r):0\le x\le1,\ \theta_r\le-x\}$ to that relaxation even yields an ideal formulation.

\subsection{Valid inequalities and network flows}\label{sec:inequalities-flows}

We now characterize all valid inequalities for $\mathcal S$. 
Consider valid inequalities for $\conv(\mathcal S)$ of the form
\begin{equation}\label{eq:network-valid-inequality}
    \bm{\beta}^{\top}\bm{\theta}\le\bm{\pi}^{\top}\bm{x}+\pi_0.
\end{equation}
With $\bm{\beta}=\bm{0}$, those are precisely the valid inequalities for $\proj_{\bm{x}}\mathcal S=\mathcal X$.
Note that we assume only that the leader has a linear formulation, thus there is no structure to be exploited in generating the valid inequalities with $\bm{\beta}=\bm{0}$, and one can resort to general integer-programming techniques.
Our focus is on inequalities with $\bm{\beta}\ne\bm{0}$, which account for the interaction between leader decisions and node potentials.
Since $\theta_t=0$, we may choose $\beta_t=-\sum_{u\in \mathcal{V}\setminus\{t\}}\beta_u$ without changing the inequality.
Throughout, we use this representative, so that $\sum_{u\in\mathcal{V}}\beta_u=0$.

All valid inequalities of interest have node potential coefficients $\bm{\beta}$ that define valid supplies over the network.
Intuitively, $\bm{\beta}$ derives from nonnegative combinations of the follower's dual inequalities, thus columns of $D$.
Combining these columns with nonnegative weights can therefore be viewed as assigning flows to the arcs.
We define the supply cone and its associated flow sets by
\begin{equation*}\label{eq:supply-cone-flows}
    \mathcal B=\{D\bm\lambda:\bm\lambda\in\R_+^{\mathcal A}\},
    \qquad
    \Lambda(\bm{\beta})=
       \{\bm{\lambda}\in\R_+^{\mathcal{A}}:D\bm{\lambda}=\bm{\beta}\}.
\end{equation*}
For $\bm{\beta}\in\mathcal B$, acyclicity makes $\Lambda(\bm{\beta})$ a nonempty polytope, and linear programming duality gives the finite, attained value
\begin{equation}\label{eq:fixed-supply-value}
    f_{\bm{\beta}}(\bm{x})
        :=\max_{\bm{\theta}\in\mathcal P(\bm{x})}
                  \bm{\beta}^{\top}\bm{\theta}
        =\min_{\bm{\lambda}\in\Lambda(\bm{\beta})}
                  \bm{q}(\bm{x})^{\top}\bm{\lambda},
        \qquad \;\forall\; \bm{x}\in\conv(\mathcal X).
\end{equation}

\begin{proposition}\label{prop:flow-validity}
    Let $\bm{\beta}$ satisfy $\sum_{u\in\mathcal V}\beta_u=0$.
    Inequality \eqref{eq:network-valid-inequality} is valid for $\mathcal S$ if and only if $\bm{\beta}\in\mathcal B$ and
    \begin{equation*}\label{eq:fixed-supply-majorant}
        \bm{\pi}^{\top}\bm{x}+\pi_0\ge f_{\bm{\beta}}(\bm{x})
        \qquad\;\forall\;\bm{x}\in\mathcal X.
    \end{equation*}
    Equivalently, for every $\bm{x}\in\mathcal X$ there exists a flow certificate $\bm{\lambda}^{\bm{x}}\in\Lambda(\bm{\beta})$ such that
    \begin{equation}\label{eq:fiber-flow-certificate}
        \bm{\pi}^{\top}\bm{x}+\pi_0
            \ge\bm{q}(\bm{x})^{\top}\bm{\lambda}^{\bm{x}}.
    \end{equation}
\end{proposition}

By \autoref{prop:flow-validity}, the right-hand side of any valid inequality defines an affine upper bound for $f_{\bm{\beta}}$ on $\mathcal X$.
Since the graph is acyclic, $\Lambda(\bm{\beta})$ is a polytope, and $f_{\bm{\beta}}$ can be expressed as a minimum of finitely many affine functions of $\bm{x}$, one for each extreme point of $\Lambda(\bm{\beta})$.
Although $f_{\bm{\beta}}$ is therefore concave, we are interested in upper bounds for the function over $\mathcal X$ only, which means they might imply stronger bounds at fractional leader coordinates than the function itself.
For $\bm{\beta}\in \mathcal B$, let $\hat{f}_{\bm{\beta}}$ be the smallest concave function on $\conv\!\left( \mathcal X \right)$ bounded below by $f_{\bm{\beta}}$ on $\mathcal X$, that is, let $\hat{f}_{\bm{\beta}}$ be the \emph{concave envelope} of $f_{\bm{\beta}}$ over $\X$.

\begin{corollary}\label{cor:fixed-supply-convexification}
    For $\bm{\beta}\in \mathcal B$, \eqref{eq:network-valid-inequality} is valid for $\mathcal S$ if and only if $\bm{\pi}^{\top}\bm{x}+\pi_0$ is an affine upper bound of $\hat f_{\bm{\beta}}$ on $\conv(\mathcal X)$.
\end{corollary}

Consequently, at every $\bar{\bm{x}}\in\conv(\mathcal X)$, there exists a supporting hyperplane of the hypograph of $\hat f_{\bm{\beta}}$ whose associated valid inequality has minimal right-hand side at $\bar{\bm{x}}$ among all valid inequalities with potential coefficients $\bm{\beta}$.

\subsection{Root supply and optimality gap}\label{sec:root-supplies}

The root potential $\theta_r$ is the only potential in the follower's optimality constraint of~\eqref{eq:single-level-reformulation}, and no potential appears in the leader's objective.
This motivates studying inequalities involving only $\bm{x}$ and $\theta_r$.
We show that this family suffices to recover the optimal value of the bilevel problem.

The \emph{root supply} $\bm{\beta}^{r}=\bm{e}_r-\bm{e}_t$ supplies one unit of flow at the root $r$ and consumes it at the target $t$.
Since $\Lambda(\bm{\beta}^r)=\conv(\mathcal Y)$, \eqref{eq:fixed-supply-value} gives
\begin{equation*}\label{eq:root-supply-value} 
    f_{\bm{\beta}^r}(\bm{x})
       =\min_{\bm{y}\in \mathcal Y}\bm{q}(\bm{x})^{\top}\bm{y}
       =\phi(\bm{x}).
\end{equation*}
Thus, the associated \emph{root cuts}\footnote{Here ``root'' refers to the root node $r$ of the graph, not to the root of a branch-and-bound tree.} have the form
$\theta_r\le\bm{\pi}^{\top}\bm{x}+\pi_0$, where the right-hand side is an affine upper bound of the follower's value function on $\mathcal X$.
Let $\hat\phi=\hat f_{\bm{\beta}^r}$ denote its concave envelope over $\X$.
By \autoref{cor:fixed-supply-convexification}, adding all root cuts to $\mathcal P$ gives the closure
\begin{equation}\label{eq:root-closure}
    \setR_r=\{(\bm{x},\bm{\theta})\in \mathcal P:\theta_r\le\hat\phi(\bm{x})\}.
\end{equation}

\begin{lemma}\label{lemma:root-projection}
    The root closure satisfies
    \begin{equation*}\label{eq:root-projection}
    \begin{aligned}
        \proj_{(\bm{x},\theta_r)}\setR_r
        &=\proj_{(\bm{x},\theta_r)}\conv(\mathcal S)=\{(\bm{x},\eta):\bm{x}\in\conv(\mathcal X),\ \eta\le\hat\phi(\bm{x})\}.
    \end{aligned}
    \end{equation*}
\end{lemma}

The projection equality implies that, for every $\bm{a}\in\R^n$ and $\gamma\le0$,
\begin{equation}\label{eq:root-objective-gap}
    \min_{(\bm{x},\bm{\theta})\in \setR_r}
        \bm{a}^{\top}\bm{x}+\gamma\theta_r
    =\min_{(\bm{x},\bm{\theta})\in \mathcal S}
        \bm{a}^{\top}\bm{x}+\gamma\theta_r.
\end{equation}
Thus, root cuts close the gap for objectives on $\mathcal S$ involving the potentials only through $\theta_r$.
This has an important consequence for problems with interdiction objectives~\cite{denegreInterdictionDiscreteBilevel2011}, which can be formulated~\cite{capraraBilevelKnapsackInterdiction2016,fischettiInterdiction2019,buiSolvingCombinatorialPricing2025} as in \eqref{eq:root-objective-gap} with $\bm{a}=0$ and $\gamma=-1$.
Precisely, for bilevel shortest-path problems with interdiction objectives, root cuts yield an LP relaxation that recovers the optimal objective value. 

More interestingly, this same family actually suffices for the general case, which we prove below.
Define the relaxation
\begin{equation*}\label{eq:root-bilevel-relaxation}
    \mathcal E_r=\{(\bm{x},\bm{\theta},\bm{y}):
        (\bm{x},\bm{\theta})\in \setR_r,\ \bm{y}\in\conv(\mathcal Y),\
        \bm{q}(\bm{x})^{\top}\bm{y}\le\theta_r\}.
\end{equation*}
For a point of $\mathcal E_r$, take a convex decomposition of its leader component attaining the root envelope and a convex decomposition of its follower path component.
The $\hat{\phi}(\bm{x})$ bound and the bilinear optimality constraint force the follower optimality gaps to zero (under the product weights), so every (positively weighted) leader--path pair is bilevel feasible.

\begin{theorem}\label{thm:root-closes-bilevel-gap}
    Every point of $\proj_{(\bm{x},\bm{y})}\mathcal E_r$ is a convex combination of bilevel-feasible leader--path pairs.
    Moreover, for every linear objective $\bm{c}^{\top}(\bm{x},\bm{y})$,
    \begin{equation*}\label{eq:root-bilevel-gap}
        \min_{(\bm{x},\bm{\theta},\bm{y})\in \mathcal E_r}
             \bm{c}^{\top}(\bm{x},\bm{y})
        =\min_{(\bm{x},\bm{\theta},\bm{y})\in\F}
             \bm{c}^{\top}(\bm{x},\bm{y}).
    \end{equation*}
\end{theorem}

The theorem strengthens the conclusion of \autoref{thm:exact-convexification} for the objective of interest: convexifying the entire leader--potential set $\mathcal S$ is unnecessary to close the bilevel gap.
The result concerns the complete family of root cuts with the bilinear optimality constraint retained; it does not assert gap closure after McCormick relaxation.

Even the complete family of root cuts can fail to describe $\conv(\mathcal S)$.
This limitation persists for the larger family of all \emph{unit supplies}, namely $\bm{e}_u-\bm{e}_v$ for distinct nodes $u,v$ joined by a directed $u\leadsto v$ path, as the following example shows.

\begin{example}[Unit supplies are insufficient]\label{eg:node-pair-cuts-insufficient}
    Let $\mathcal X=\{0,1\}^2$ and $\mathcal{V}=\{r,u,v,t\}$.
    The arcs $r\to u$ and $r\to v$ have cost zero.
    There are two parallel $u\to t$ arcs of costs $0$ and $1-x_1-x_2$, and two parallel $v\to t$ arcs of costs $0$ and $x_2-x_1$, as illustrated in Figure~\ref{fig:example-networks}.
    The point
    \[
        \bar{\bm{x}}=(1/2,1/2),\qquad
        (\bar\theta_r,\bar\theta_u,\bar\theta_v,\bar\theta_t)
            =(-1/2,0,0,0)
    \]
    belongs to $\mathcal P$ and satisfies every cut with a unit supply, including every root cut.
    Nevertheless, the inequality
    \begin{equation*}\label{eq:coupled-frontier-cut}
        \theta_u+\theta_v\le -x_1
    \end{equation*}
    is valid for $\mathcal S$ and violated by this point.
    Its supply is $\bm{e}_u+\bm{e}_v-2\bm{e}_t$.
\end{example}

Unit supplies contain a generator of every extreme ray of $\mathcal B$, since this cone is generated by the columns of $D$, each of which is a unit supply.
Thus even inequalities from all extreme supply directions can leave points outside $\conv(\mathcal S)$ unseparated.

\section{Cut Generation}\label{sec:cut-generation}

In this section, we explore the generation of linear valid inequalities (cuts) to separate a given point of a continuous relaxation of~\eqref{eq:single-level-reformulation}, possibly already strengthened by
other valid inequalities.
In particular, we study cut generation for the leader--potential set, construct node potential coefficients that separate fractional vertices of \(\mathcal P\), and analyze the complexity of separation under partially-fixed coefficients.

Our structured approach to the single-level reformulation \eqref{eq:single-level-reformulation} can guide the choice of variables involved in separation.
Let $(\bar{\bm{x}},\bar{\bm{\theta}},\bar{\bm{y}})$ be an optimal solution from any continuous relaxation of $\F$.
If $\bar{\bm{x}}$ is fractional and its active leader constraints have rank $n$, then an inequality for $\mathcal X$ alone suffices to separate it.
As previously discussed, we will not concentrate on this case.

For candidates with $\overline{\bm{x}}\in \conv(\mathcal X)$, \autoref{thm:exact-convexification} justifies focusing on the leader--potential set.
Define $\Xi$ to be a continuous relaxation of $\Omega$ as follows:
\[\Xi:= \left\{(\bm{x}, \bm{\theta}, \bm{y})\in \conv(\mathcal X)\times\R^{\mathcal{V}} \times \conv(\mathcal Y): \; \bm{q}(\bm{x})^{\top}\bm{y}\le\theta_r, \; D^{\top}\bm{\theta}\le\bm{q}(\bm{x}),\; \theta_t=0\right\}.\] 
The proposition also shows that inequalities of the form \eqref{eq:network-valid-inequality} with $\bm{\beta}\ne\bm{0}$ are sufficient to separate any point in $\Xi \setminus \conv\!\left( \F  \right)$.
We explore this result in the following lemma.

\begin{lemma}\label{lemma:bilinear-extreme-points}
    Every extreme point $(\bar{\bm{x}},\bar{\bm{\theta}},\bar{\bm{y}})\in\Xi$ satisfies $\bar{\bm{y}}\in \mathcal Y$, and $(\bar{\bm{x}},\bar{\bm{\theta}})$ is a vertex of $\mathcal P$.
\end{lemma}

We therefore focus on vertices of \(\mathcal P\) with fractional leader coordinates, which necessarily lie outside \(\conv(\mathcal S)\) and admit a separating inequality.

\subsection{Fixing node potential coefficients $\bm \beta$}\label{sec:supply-selection} 

Let $(\bar{\bm{x}},\bar{\bm{\theta}})$ be a vertex of $\mathcal P$ with fractional leader coordinates.
By \autoref{prop:flow-validity}, a cut-generating linear program (CGLP) can jointly select $(\bm{\beta},\bm{\pi},\pi_0)$ and the flow certificates, subject to a linear normalization fixing the coefficient scale.
Its explicit formulation requires one flow certificate \eqref{eq:fiber-flow-certificate} per leader assignment, hence $|\X||\mathcal A|$ flow variables, and may be exponentially large.
We instead construct the node potential coefficients $\bm \beta$ from the arc inequalities active at $(\bar{\bm{x}},\bar{\bm{\theta}})$, preserving the existence of a separating cut while leaving the leader coefficients to be optimized.

\autoref{prop:flow-validity} restricts admissible $\bm{\beta}$ to the supply cone $\mathcal B$, but does not determine which supplies admit separation at the candidate.
For fixed $\bm{\beta}\in \mathcal B\setminus\{\bm{0}\}$, \autoref{cor:fixed-supply-convexification} gives the maximum attainable violation as
\begin{equation}\label{eq:fixed-supply-violation}
    v_{\bm{\beta}}(\bar{\bm{x}},\bar{\bm{\theta}})
      :=\bm{\beta}^{\top}\bar{\bm{\theta}}
         -\hat f_{\bm{\beta}}(\bar{\bm{x}}).
\end{equation}
A separating inequality with this $\bm\beta$ exists if and only if this value is positive.

For any $(\bar{\bm{x}},\bar{\bm{\theta}})\in \mathcal P$, define its active arcs by
\begin{equation*}\label{eq:active-arcs}
    \mathcal A_0=\{a\in\mathcal A:
        D_a^{\top}\bar{\bm{\theta}}=q_a(\bar{\bm{x}})\}.
\end{equation*}
Connectivity and spanning trees below refer to the underlying undirected graph; each selected arc retains its original orientation.

\begin{lemma}\label{lemma:active-arc-connectivity}
    For $(\bar{\bm{x}},\bar{\bm{\theta}})\in \mathcal P$, the active-arc graph $(\mathcal V,\mathcal A_0)$ is connected if and only if $\bar{\bm{\theta}}$ is a vertex of $\mathcal P(\bar{\bm{x}})$.
    In particular, it is connected at every vertex of $\mathcal P$.
\end{lemma}

\begin{figure}[t]
    \centering
    \begin{tikzpicture}[scale=.9,transform shape,font=\small]
        \input{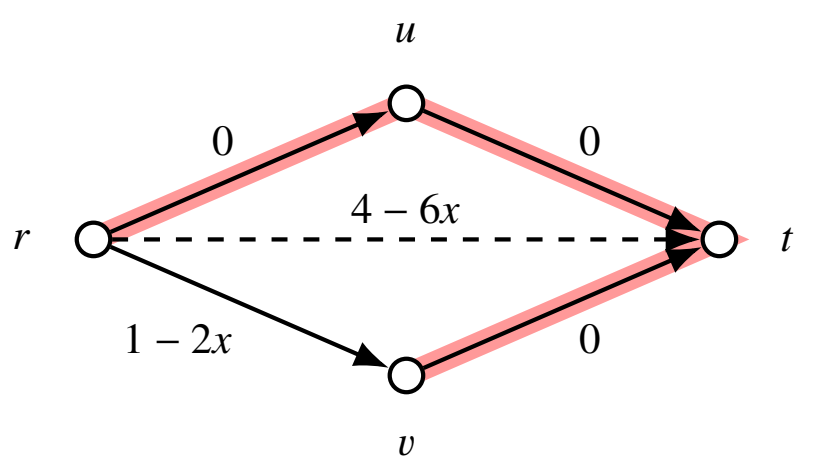}
    \end{tikzpicture}
    \caption{Active arcs at a vertex of $\mathcal P$ with $\mathcal X=\{0,1\}$ and $(\bar x,\bar{\bm\theta})=(1/2,\bm0)$. Arc labels give costs. Solid arcs are active, and the dashed arc has slack $1$. The red highlight marks a spanning tree of active arcs. The four active equalities uniquely determine $(x,\theta_r,\theta_u,\theta_v)$.}
    \label{fig:active-arcs}
\end{figure}

\autoref{fig:active-arcs} illustrates the connected active-arc graph at a fractional vertex of $\mathcal P$.

Our construction is based on a standard geometric observation: a positive combination of full-rank active inequality rows yields an optimization direction uniquely maximized by the associated basic solution.
We apply this observation to the fiber $\mathcal{P}(\bm x)$.
The network structure makes the required row selection inexpensive: after eliminating $\theta_t$, a spanning tree of active arcs provides a minimum-cardinality full-rank collection of $|\mathcal V|-1$ rows.
Once the active-arc graph is available, such a tree can be found by graph traversal in $O(|\mathcal V|+|\mathcal A_0|)$ time, without matrix rank computations.
Positive weights on its arcs then determine a supply vector $\bm\beta$; the following lemma establishes the resulting separation guarantee.

\begin{lemma}\label{lemma:active-tree-supply}
    Suppose $(\bar{\bm{x}},\bar{\bm{\theta}})\in \mathcal P$ has a connected active-arc graph. Choose a spanning tree $T\subseteq\mathcal A_0$ and weights $w_a>0$ for $a\in T$.
    Then 
    \begin{equation*}\label{eq:active-tree-supply}
        \bm{\beta}=\sum_{a\in T}w_aD_a
    \end{equation*}
    belongs to $\mathcal B\setminus\{\bm{0}\}$. Moreover, $v_{\bm{\beta}}(\bar{\bm{x}},\bar{\bm{\theta}})>0$ if and only if $(\bar{\bm{x}},\bar{\bm{\theta}})\notin\conv(\mathcal S)$.
    In particular, this construction guarantees separation at every vertex of $\mathcal P$ with fractional leader coordinates.
\end{lemma}

The weighted sum of the tree inequalities is valid for $\mathcal P$ and active at $(\bar{\bm{x}},\bar{\bm{\theta}})$, so it does not itself separate.
It selects the potential coefficients; separation comes from strengthening the affine right-hand side over the discrete leader set $\mathcal{X}$.
The guarantee concerns the existence of a cut, not its strength relative to other supply vectors.

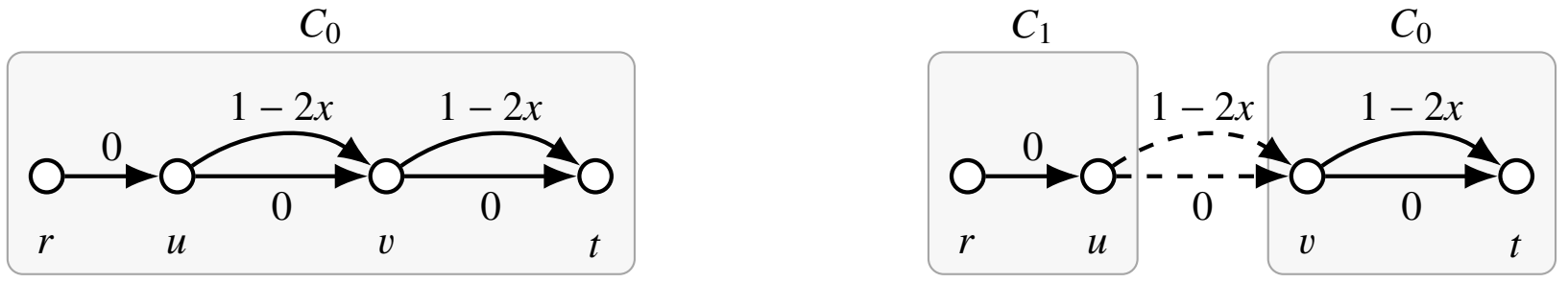
\begin{figure}[htbp]
    \centering
    \subfloat[Vertex of $\mathcal P$\label{fig:components-before}]{%
        \begin{minipage}{.48\linewidth}
            \centering
            \begin{tikzpicture}[scale=.9,transform shape,font=\small]
                \input{figures/components-before.pic.tex}
            \end{tikzpicture}
            \par\smallskip
            {\small $\bar x=1/2,\quad\bar{\bm\theta}=(0,0,0,0)$}
        \end{minipage}%
    }\hfill
    \subfloat[Vertex after adding $\theta_r\le-2x$\label{fig:components-after}]{%
        \begin{minipage}{.48\linewidth}
            \centering
            \begin{tikzpicture}[scale=.9,transform shape,font=\small]
                \input{figures/components-after.pic.tex}
            \end{tikzpicture}
            \par\smallskip
            {\small $\bar x=1/2,\quad\bar{\bm\theta}=(-1,-1,0,0)$}
        \end{minipage}%
    }
    \caption{Connected components before and after adding one cut, with $\mathcal X=\{0,1\}$ and potentials ordered as $(r,u,v,t)$. Solid arcs are active, dashed arcs are inactive, and shading identifies components of the original active-arc graph. The root cut is valid since $\phi(0)=0$ and $\phi(1)=-2$. In (b), it is active and fixes the relative shift between $C_1$ and $C_0$.}
    \label{fig:active-components}
\end{figure}

Let $\mathcal P^+$ be obtained from $\mathcal P$ by adding finitely many valid inequalities $(\bm\beta^j)^\top\bm\theta\le\ell_j(\bm x)$, with affine $\ell_j$ and $\bm\beta^j\in \mathcal B$.
Unlike for $\mathcal P$, the active-arc graph at a vertex of $\mathcal P^+$ need not be connected.
Consider first the case in which $\left( \mathcal{V},\mathcal{A}_0 \right) $ consists of two components, $C_0$ containing $t$ and $C_1$.
A spanning tree in each component fixes all differences of node potential variables, allowing only for a common shift of the potentials in $C_1$.
An active added cut fixes this shift whenever its total supply $\sum_{u\in C_1}\beta^j_u$ is nonzero.
At a vertex, at least one such cut must exist; otherwise, a small shift in either direction would preserve feasibility.
Thus two trees and one active added cut suffice to determine the node potentials.
\autoref{fig:active-components} illustrates how adding a cut can produce a new fractional vertex with two active-arc components.
The following proposition extends this construction to any number of components.

\begin{proposition}\label{prop:strengthened-supply}
    Let $(\bar{\bm x},\bar{\bm\theta})$ be a vertex of $\mathcal P^+$, let $C_0,\ldots,C_{c-1}$ be its active-arc components, with $t\in C_0$, and let $\mathcal J_0$ index its active added cuts.
    The matrix
    \[
        \Gamma_{j\ell}=\sum_{u\in C_\ell}\beta^j_u,
        \qquad j\in\mathcal J_0,\quad \ell=1,\ldots,c-1,
    \]
    has rank $c-1$.
    Choose a spanning tree in each component, with union $F$, and a set $J\subseteq\mathcal J_0$ indexing $c-1$ independent rows of $\Gamma$.
    The corresponding $|\mathcal V|-1$ original arc inequalities and added cuts form a minimum-cardinality full-rank collection in the node potential coordinates, excluding $\theta_t$.
    For any positive weights $w_a$ and $\mu_j$, the supply
    \[
        \bm\beta=\sum_{a\in F}w_aD_a+\sum_{j\in J}\mu_j\bm\beta^j
    \]
    belongs to $\mathcal B\setminus\{\bm0\}$. Moreover, $v_{\bm{\beta}}(\bar{\bm{x}},\bar{\bm{\theta}})>0$ if and only if $(\bar{\bm{x}},\bar{\bm{\theta}})\notin\conv(\mathcal S)$.
    In particular, it admits a separating inequality whenever $\bar{\bm x}$ is fractional.
\end{proposition}

The network structure therefore keeps row selection inexpensive within each component: graph traversal selects the forest, and $c-1$ linearly independent active cuts fix the remaining degrees of freedom.
For $c=1$, no added cut is needed and the spanning-tree construction applies directly.

On a practical setting, there is strong reason to start with root supplies before attempting the constructions suggested above.
By \autoref{thm:root-closes-bilevel-gap}, root cuts suffice to close the bilevel objective gap when the bilinear optimality constraint is retained.
Moreover, each root cut has the node potential coefficients of an additional \(r\to t\) arc, so adding root cuts preserves the spanning-tree construction on the augmented active graph (\autoref{lemma:active-tree-supply}).
Root cuts can therefore be pursued first while retaining an inexpensive way to select other supply vectors when needed.

\subsection{Complexity of separation}\label{sec:separation-complexity}

Selecting a supply vector $\bm\beta$ that admits separation leaves the leader variable coefficients to be determined.
For fixed $\bm{\beta}\in \mathcal B\setminus\{\bm0\}$ and a point $(\bar{\bm{x}},\bar{\bm{\theta}})\in \mathcal P$, \autoref{cor:fixed-supply-convexification} gives the attained maximum violation as
\begin{equation}\label{eq:fixed-supply-separation}
\begin{aligned}
    v_{\bm{\beta}}(\bar{\bm{x}},\bar{\bm{\theta}})
      =\max_{\bm{\pi},\pi_0}\quad
        &\bm{\beta}^{\top}\bar{\bm{\theta}}
           -\bm{\pi}^{\top}\bar{\bm{x}}-\pi_0\\
    \textrm{s.t.}\quad
        &\bm{\pi}^{\top}\bm{x}^k+\pi_0
           \ge f_{\bm{\beta}}(\bm{x}^k),
           \qquad k=1,\ldots,K,
\end{aligned}
\end{equation}
where $\bm x^1,\ldots,\bm x^K$ is an enumeration of the feasible points in $\mathcal{X}$.
Each value $f_{\bm{\beta}}(\bm{x}^k)$ can be computed by a minimum-cost flow problem.
Given these values, LP~\eqref{eq:fixed-supply-separation} has only $n+1$ variables, but still one constraint per leader assignment.
Replacing each constraint by its flow certificate \eqref{eq:fiber-flow-certificate} gives an equivalent CGLP without precomputing the values.
The resulting formulation can however still be exponential in size.

We propose two ways to address the exponential size of the CGLP.
With respect to the number of leader assignments, we can, for example, consider a split disjunction on one of the fractional coordinates~\cite{Balas2018LnP}.
With respect to the size of the underlying graph, we can restrict the minimum-cost flow problem to the active arcs $\mathcal A_0$.
For a vertex of $\mathcal{P}$ with fractional leader coordinates and a supply constructed by \autoref{lemma:active-tree-supply}, these reductions also preserve the existence of a separating inequality.
The two reductions can decrease the maximum violation compared with the original CGLP with fixed $\bm\beta$.
When a linear description of $\conv(\mathcal X)$ is supplied, the reduced CGLP has polynomial size in that description and the active graph.

To motivate the split relaxation, we also examine the computational complexity of the separation problem for a particular fixed supply vector on a small graph.
By \emph{exact separation}, we mean to determine occurrence of a violation, and if so, to find a violated cut with the fixed supply vector $\bm{\beta}$.

\begin{proposition}\label{prop:separation-hardness}
    Membership in the root-cut closure $\setR_r$ is NP-hard, even for candidates in $\mathcal P$ when $\mathcal X=\{0,1\}^n$ and the graph consists of two parallel $r\to t$ arcs, one with constant cost.
    Consequently, exact separation of root cuts and exact evaluation of $\hat \phi$ are NP-hard under the same restrictions.
\end{proposition}

The proof, given in \autoref{appx:proofs}, reduces \textsc{Partition} to membership at $\bar{\bm{x}}=\tfrac12\bm1$.
Both arc inequalities are active at the constructed candidate, so restricting to active arcs leaves the hardness construction intact.
The result concerns arbitrary candidates in $\mathcal P$ and also establishes hardness of general fixed-supply separation.

Exact separation therefore remains NP-hard even with an unconstrained leader and a graph consisting of two parallel arcs.
The split relaxation addresses a computational difficulty that persists after graph reduction.
Combined with the supply construction of \autoref{lemma:active-tree-supply}, it retains a separating inequality at fractional vertices while admitting a CGLP of polynomial size in the supplied leader formulation and the active graph. To see this, each branch of the split disjunction is described by the given linear description of $\conv(\X)$, the retained arc inequalities, and one split equation. Dualizing cut validity on each disjunction branch introduces one multiplier per row, one coefficient-matching equation per original variable, and one inequality for the right-hand side. Since the disjunction has two branches, the resulting CGLP has size polynomial in the supplied leader description and the active graph.

\section{Conclusions}

In this work, motivated by network pricing, follower's dynamic programming formulations and certain interdiction models, we developed a polyhedral characterization of the bilevel 0--1 shortest-path problem by studying the interaction between leader decisions and follower node potentials.
Network flows characterize the valid inequalities for this joint set.
With the leader set convexified and the bilinear optimality constraint retained, root cuts strengthen the continuous relaxation of the problem and recover the bilevel optimum for every linear leader objective.
Full convexification of the leader--potential set can require coupled supplies even after all unit-supply cuts are added.

We proposed cut generation by constructively fixing supplies.
Spanning trees of active arcs yield separating supplies at vertices of $\mathcal P$ with fractional leader coordinates, and active-arc components extend this construction after cuts are added.
Nonetheless, we showed that exact fixed-supply separation is NP-hard even for networks with two parallel arcs.
For the initial relaxation, split disjunctions and restriction to active arcs reduce the CGLP while preserving separation with the constructed supplies.

\section*{Acknowledgments}

This work was funded by the FRQNT Doctoral Scholarship No. B2X-367536, and the NSERC Grants No. 2024-04051 and 2024-05059.

OpenAI Codex (GPT-6) was used to assist with literature searches, \LaTeX{} edits, and figure generation.
The authors reviewed the content as needed and take full responsibility for this manuscript.

\bibliographystyle{splncs04}
\bibliography{references}

\appendix

\section{Deferred Proofs}\label{appx:proofs}

\subsection{Polyhedrality of the single-level feasible hull}\label{appx:proof-polyhedrality}

\begin{proofof}{\autoref{prop:SP-is-polyhedral}}
    Fix $(\hat{\bm{x}},\hat{\bm{y}})\in\mathcal X\times\mathcal Y$ for which the potential constraints in \eqref{eq:single-level-reformulation} are feasible.
    Weak duality and the optimality constraint give
    $\theta_r=\bm{q}(\hat{\bm{x}})^{\top}\hat{\bm{y}}$.
    For each node $u$, choose an $r\leadsto u$ path $\rho_u$ and a $u\leadsto t$ path $\tau_u$, allowing the empty path at an endpoint.
    Summing the dual arc inequalities along these paths yields
    \[
        \bm{q}(\hat{\bm{x}})^{\top}\hat{\bm{y}}
            -\sum_{a\in\rho_u}q_a(\hat{\bm{x}})
        \le\theta_u\le\sum_{a\in\tau_u}q_a(\hat{\bm{x}}).
    \]
    Thus the feasible potentials for each fixed pair form a bounded polyhedron.
    Since $\mathcal X$ and $\mathcal Y$ are finite, $\F$ is a finite union of polytopes, and its convex hull is a polytope.
\end{proofof}

\subsection{Effect of convexifying the dual network polyhedron}\label{appx:proof-exact-convexification}

\begin{proofof}{\autoref{thm:exact-convexification}}
    It suffices to show $\F \subseteq \mathcal E \subseteq \conv(\F)$.
    Taking the convex hull gives $\conv(\mathcal E)=\conv(\F)$.
    
    The inclusion $\F\subseteq \mathcal E$ follows from $\mathcal S\subseteq\conv(\mathcal S)$ and $\mathcal Y\subseteq\conv(\mathcal Y)$.
    For the reverse containment in $\conv(\F)$, take $(\bar{\bm{x}},\bar{\bm{\theta}},\bar{\bm{y}})\in \mathcal E$ and choose finite convex representations
    \[
        (\bar{\bm{x}},\bar{\bm{\theta}})
          =\sum_i\lambda_i(\bm{x}^i,\bm{\theta}^i),
        \qquad
        \bar{\bm{y}}=\sum_j\mu_j\bm{y}^j,
    \]
    where $(\bm{x}^i,\bm{\theta}^i)\in \mathcal S$, $\bm{y}^j\in \mathcal Y$, $\lambda_i,\mu_j\ge0$, and $\sum_i\lambda_i=\sum_j\mu_j=1$.
    Since the path set $\mathcal Y$ is independent of the leader assignment, every pair $(i,j)$ combines a feasible follower path with dual-feasible potentials at $\bm{x}^i$.
    Weak duality then gives us $\bm{q}(\bm{x}^i)^{\top}\bm{y}^j-\theta_r^i\ge 0$ for each pair $(i,j)$.
    Because the candidate solution satisfies the follower's optimality constraint and $\bm{q}(\bm{x})$ is affine, we have
    \begin{equation*}
        0\ge\bm{q}(\bar{\bm{x}})^{\top}\bar{\bm{y}}-\bar\theta_r = \sum_{i,j} \lambda_i \mu_j (\bm{q}(\bm{x}^{i})^{\top}\bm{y}^{j} - \theta_r^{i}) \ge 0
    .\end{equation*}
    For every $(i,j)$ pair with $\lambda_i \mu_j >0$ the optimality gap must vanish, that is, $\bm{q}(\bm{x}^i)^{\top}\bm{y}^j-\theta_r^i= 0$.
    Hence each positively weighted triple $(\bm{x}^i,\bm{\theta}^i,\bm{y}^j)$ belongs to $\F$.
    Finally, the weights $\lambda_i\mu_j$ sum to one and satisfy
    \[
        (\bar{\bm{x}},\bar{\bm{\theta}},\bar{\bm{y}})
        =\sum_{i,j}\lambda_i\mu_j
          (\bm{x}^i,\bm{\theta}^i,\bm{y}^j),
    \]
    proving $\mathcal E\subseteq\conv(\F)$.
\end{proofof}

\subsection{Flow characterization and fixed-supply convexification}\label{appx:proof-flow-characterization}

\begin{proofof}{\autoref{prop:flow-validity}}
    We first verify \eqref{eq:fixed-supply-value} and supply admissibility.
    For any $\bm{x}\in\conv(\mathcal X)$, the potential maximization is feasible, and its linear programming dual is
    \[
        \min\{\bm{q}(\bm{x})^{\top}\bm{\lambda}:
              D\bm{\lambda}+\mu\bm{e}_t=\bm{\beta},\
              \bm{\lambda}\ge\bm0,\ \mu\in\R\}.
    \]
    Summing the balance equations gives $\mu=0$ because $\sum_{u\in \mathcal{V}}\beta_u=0$.
    Hence dual feasibility is equivalent to $\bm{\beta}\in\mathcal B$.
    If \eqref{eq:network-valid-inequality} is valid, its right-hand side bounds the potential objective on any fiber indexed by $\bm{x}\in\mathcal X$.
    Such a fiber exists because $\mathcal X\ne\varnothing$, so linear programming duality implies $\bm{\beta}\in\mathcal B$.

    For $\bm{\beta}\in\mathcal B$, the recession cone of $\Lambda(\bm{\beta})$ consists of nonnegative circulations.
    Acyclicity implies that the only such circulation is zero, so $\Lambda(\bm{\beta})$ is a nonempty polytope.
    Strong duality now gives \eqref{eq:fixed-supply-value}, with both optima finite and attained.
    Maximizing the left-hand side of \eqref{eq:network-valid-inequality} on each integer fiber proves the affine-upper-bound characterization.
    An optimal flow on that fiber supplies \eqref{eq:fiber-flow-certificate} whenever the upper bound condition holds.
    Conversely, any such certificate gives
    \[
        \bm{\beta}^{\top}\bm{\theta}
        =(\bm{\lambda}^{\bm{x}})^{\top}D^{\top}\bm{\theta}
        \le\bm{q}(\bm{x})^{\top}\bm{\lambda}^{\bm{x}}
        \le\bm{\pi}^{\top}\bm{x}+\pi_0
        \qquad\text{for all }\bm{\theta}\in\mathcal S(\bm{x}),
    \]
    proving validity.
\end{proofof}

\subsection{Root supply and the optimality gap}\label{appx:proof-root-supplies}

\begin{proofof}{\autoref{lemma:root-projection}}
    We first determine the possible values of the root potential $\theta_r$ for a fixed leader decision $\bm{x}\in\conv(\mathcal X)$.
    Assigning each node its shortest-path cost to $t$ gives feasible potentials with $\theta_r=\phi(\bm{x})$.
    Since $r$ has no incoming arcs, decreasing only $\theta_r$ preserves feasibility.
    Together with weak duality, this proves
    \[
        \proj_{\theta_r}\mathcal P(\bm{x})=(-\infty,\phi(\bm{x})].
    \]

    We next show that convexifying $\mathcal S$ replaces this upper bound by the concave envelope of the values at integer leader decisions.
    Let $\bm x^1,\ldots,\bm x^K$ be an enumeration of $\mathcal{X}$ and fix $\bar{\bm x}\in\conv(\mathcal X)$.
    The concave envelope has the representation
    \begin{equation}\label{eq:app-envelope-lp}
        \hat\phi(\bar{\bm x})
        =\max\left\{\sum_{k=1}^K\alpha_k\phi(\bm x^k):
            \sum_{k=1}^K\alpha_k\bm x^k=\bar{\bm x},\
            \sum_{k=1}^K\alpha_k=1,\ \bm\alpha\ge\bm0\right\}.
    \end{equation}
    Indeed, the right-hand side is a concave upper bound of the values on $\mathcal X$, and every such upper bound dominates it by concavity.
    The maximum is attained because the feasible weight set is nonempty and compact and each $\phi(\bm x^k)$ is finite.
    Every convex combination of points of $\mathcal S$ with leader component $\bar{\bm x}$ has root potential at most the corresponding weighted sum of follower values, and hence at most $\hat\phi(\bar{\bm x})$.
    To attain this bound, choose weights $\alpha_k$ maximizing \eqref{eq:app-envelope-lp} and, for each $k$, let $\bm\theta^k$ be the shortest-path potentials for arc costs $\bm q(\bm x^k)$.
    Then $(\bm x^k,\bm\theta^k)\in\mathcal S$ and $\theta_r^k=\phi(\bm x^k)$, so
    \[
        (\bar{\bm x},\bar{\bm\theta})
            :=\sum_k\alpha_k(\bm x^k,\bm\theta^k)\in\conv(\mathcal S),
        \qquad \bar\theta_r=\hat\phi(\bar{\bm x}).
    \]
    Decreasing every $\theta_r^k$ by the same nonnegative amount preserves membership in $\mathcal S$ and yields every smaller value of $\bar\theta_r$.
    This proves the asserted projection of $\conv(\mathcal S)$.

    Finally, we compare this projection with that of the root closure.
    Since $\phi$ is concave, $\hat\phi\le\phi$ on $\conv(\mathcal X)$.
    The first part of the proof therefore shows that imposing $\theta_r\le\hat\phi(\bm x)$ on $\mathcal P$ leaves precisely the root potential values $(-\infty,\hat\phi(\bm x)]$ for each $\bm x\in\conv(\mathcal X)$.
    By \eqref{eq:root-closure}, this proves the equality for $\setR_r$.
\end{proofof}

\begin{proofof}{\autoref{thm:root-closes-bilevel-gap}}
    Take $(\bar{\bm{x}},\bar{\bm{\theta}},\bar{\bm{y}})\in\mathcal E_r$.
    Choose weights $\alpha_k$ attaining \eqref{eq:app-envelope-lp} and write $\bar{\bm{y}}=\sum_j\mu_j\bm{y}^j$ with $\bm{y}^j\in\mathcal Y$, $\mu_j\ge0$, and $\sum_j\mu_j=1$.
    Feasibility and the affinity of $\bm q$ give
    \[
    \begin{aligned}
        \hat\phi(\bar{\bm{x}})
        &\ge\bar\theta_r
         \ge\bm q(\bar{\bm{x}})^{\top}\bar{\bm{y}}\\
        &=\sum_{k,j}\alpha_k\mu_j
                        \bm q(\bm{x}^k)^{\top}\bm{y}^j
         \ge\sum_k\alpha_k\phi(\bm{x}^k)
         =\hat\phi(\bar{\bm{x}}).
    \end{aligned}
    \]
    All inequalities are equalities.
    Each difference $\bm q(\bm{x}^k)^{\top}\bm{y}^j-\phi(\bm{x}^k)$ is nonnegative, so it vanishes whenever $\alpha_k\mu_j>0$.
    Every positively weighted pair $(\bm{x}^k,\bm{y}^j)$ is therefore bilevel feasible, and
    \[
        (\bar{\bm{x}},\bar{\bm{y}})
        =\sum_{k,j}\alpha_k\mu_j(\bm{x}^k,\bm{y}^j).
    \]
    This proves the convex-combination claim and bounds every linear leader objective below by the bilevel optimum.
    Conversely, $\F\subseteq\mathcal E_r$ because root cuts are valid for $\mathcal S$.
    The optimal values consequently coincide.
\end{proofof}

\subsection{Reduction to the dual network polyhedron}\label{appx:proof-separation-reduction}

\begin{proofof}{\autoref{lemma:bilinear-extreme-points}}
    Let $(\bar{\bm{x}},\bar{\bm{\theta}},\bar{\bm{y}})$ be an extreme point of $\mathcal R$.
    Decompose $\bar{\bm{y}}$ into path-incidence vectors.
    Weak duality gives $\bm q(\bar{\bm{x}})^{\top}\bm y\ge\bar\theta_r$ for each path, and the optimality constraint forces equality in their weighted average.
    Every positively weighted path therefore gives a feasible triple with the same $(\bar{\bm{x}},\bar{\bm{\theta}})$.
    Extremality forces all those paths to equal $\bar{\bm{y}}$, proving $\bar{\bm{y}}\in\mathcal Y$.

    Now suppose $(\bar{\bm{x}},\bar{\bm{\theta}})$ is a nontrivial convex combination of two points of $\mathcal P$.
    With $\bar{\bm{y}}$ fixed, the gap $\bm q(\bm{x})^{\top}\bar{\bm{y}}-\theta_r$ is affine and nonnegative on $\mathcal P$.
    It vanishes at the candidate, so it vanishes at both endpoints.
    Appending $\bar{\bm{y}}$ to these endpoints gives a nontrivial convex representation in $\mathcal R$, a contradiction.
\end{proofof}

\subsection{Constructive supply selection}\label{appx:proof-supply-selection}

\begin{proofof}{\autoref{lemma:active-arc-connectivity}}
    If the active-arc graph of a point $(\bar{\bm x},\bar{\bm \theta})\in\mathcal{P}$ is connected, a displacement of node potentials preserving its arc equalities must be constant on all nodes.
    The normalization at $t$ makes that displacement zero.
    Hence the active equalities uniquely determine the node potentials from $\bar{\bm x}$, proving that $\bar{\bm{\theta}}$ is a vertex of the fiber $\mathcal{P}(\bar{\bm x})$.
    If the graph is disconnected, shift every node potential in a component not containing $t$ by $\epsilon$ or by $-\epsilon$.
    Active equalities are preserved, and the finitely many inactive arc inequalities remain feasible for sufficiently small $\epsilon>0$.
    These two distinct feasible node potential vectors have midpoint $\bar{\bm{\theta}}$, contradicting fiber extremality.
    Finally, a nontrivial convex representation within a fiber also gives one in $\mathcal P$, so every vertex of $\mathcal P$ is a fiber vertex.
\end{proofof}

\begin{proofof}{\autoref{lemma:active-tree-supply}}
    The vector $\bm{\beta}$ is the supply of the nonnegative flow with value $w_a$ on $T$ and zero elsewhere.
    It is nonzero because a DAG has no nonzero nonnegative circulation.
    Summing the selected inequalities gives the valid inequality
    \[
        \bm{\beta}^{\top}\bm{\theta}\le\ell(\bm{x}),
        \qquad \ell(\bm{x})=\sum_{a\in T}w_aq_a(\bm{x}),
    \]
    with equality at the candidate.

    The node potentials $\bar{\bm{\theta}}$ uniquely maximize $\bm{\beta}^{\top}\bm{\theta}$ on $\mathcal P(\bar{\bm{x}})$.
    To see this, fix $\bm{x}=\bar{\bm{x}}$ and take any $\bm{\theta}\in\mathcal P(\bar{\bm{x}})$.
    The slack of the aggregated inequality is
    \[
        \ell(\bar{\bm{x}})-\bm{\beta}^{\top}\bm{\theta}
        =\sum_{a\in T}w_a\bigl(q_a(\bar{\bm{x}})-D_a^{\top}\bm{\theta}\bigr)
        \ge0.
    \]
    The candidate $\bar{\bm{\theta}}$ attains the upper bound $\ell(\bar{\bm{x}})$ because every arc of $T$ is active there.
    Any potential vector $\bm{\theta}$ attaining the same bound must make every slack in the sum zero, since the slacks are nonnegative and all weights are strictly positive.
    Hence $D_a^{\top}(\bm{\theta}-\bar{\bm{\theta}})=0$ for every $a\in T$, so the displacement $\bm{\theta}-\bar{\bm{\theta}}$ has equal values at the endpoints of every tree arc.
    Connectivity makes this displacement constant on all nodes, and $\theta_t=\bar\theta_t=0$ makes the constant zero.
    In other words, $\bm{\theta}$ must be equal to $\overline{\bm{\theta}}$.

    Applying the convex-combination representation in \eqref{eq:app-envelope-lp} to $f_{\bm\beta}$ gives the attained identity
    \begin{equation}\label{eq:app-hull-fiber-value}
        \hat f_{\bm\beta}(\bar{\bm x})
        =\max\{\bm\beta^{\top}\bm\theta:
                     (\bar{\bm x},\bm\theta)\in\conv(\mathcal S)\}.
    \end{equation}
    Indeed, any convex combination of points of $\mathcal S$ has objective at most the corresponding weighted integer-fiber values, and maximizing weights together with maximizing fiber potentials attain equality.
    Since $\conv(\mathcal S)\subseteq\mathcal P$, \eqref{eq:app-hull-fiber-value} gives
    $\hat f_{\bm{\beta}}(\bar{\bm{x}})\le\bm{\beta}^{\top}\bar{\bm{\theta}}$.
    Equality holds if the candidate belongs to $\conv(\mathcal S)$.
    Conversely, if equality holds, an attaining hull potential also maximizes the objective $\bm{\beta}^{\top}\bm{\theta}$ over $\mathcal P(\bar{\bm{x}})$, so uniqueness makes it $\bar{\bm{\theta}}$.
    
    The difference $\bm{\beta}^{\top}\bar{\bm{\theta}}-\hat f_{\bm{\beta}}(\bar{\bm{x}})$ is the maximum attainable violation $v_{\bm{\beta}}(\bar{\bm{x}},\bar{\bm{\theta}})$ defined in \eqref{eq:fixed-supply-violation}.
    We have shown that this difference is nonnegative and vanishes precisely when the candidate belongs to $\conv(\mathcal S)$.
    A vertex of $\mathcal P$ in $\conv(\mathcal S)$ must itself belong to $\mathcal S$, since all positively weighted points in a convex representation from $\mathcal S\subseteq\mathcal P$ must equal that vertex.
    A vertex with fractional leader coordinates is therefore outside $\conv(\mathcal S)$, proving the final assertion.
\end{proofof}

\begin{proofof}{\autoref{prop:strengthened-supply}}
    A displacement preserving all active arc equalities is constant on each component $C_\ell$ and zero on $C_0$ after imposing $\theta_t=0$.
    Write its remaining component values as $\delta_1,\ldots,\delta_{c-1}$.
    It preserves the active added cuts precisely when $\Gamma\bm\delta=\bm0$.
    If $\Gamma$ had rank less than $c-1$, a nonzero such displacement of sufficiently small magnitude in either direction would preserve every constraint of $\mathcal P^+$, contradicting the vertex assumption.
    Thus $\rank\Gamma=c-1$.

    The selected forest has $|\mathcal V|-c$ arcs and forces componentwise constancy.
    The selected $c-1$ rows of $\Gamma$ force all remaining component shifts to vanish.
    Consequently, the $|\mathcal V|-1$ selected rows have full rank in the free potential coordinates, which is the minimum possible cardinality.
    This also covers $c=1$, when no added row is needed.

    Each selected potential coefficient belongs to $\mathcal B$, so their positive combination $\bm\beta$ also belongs to $\mathcal B$.
    Choose a nonnegative flow representing each selected cut coefficient and aggregate these flows with the selected arc flows using the prescribed positive weights.
    Full rank ensures that at least one selected coefficient is nonzero, so the aggregate flow is nonzero.
    Acyclicity precludes it from being a circulation, proving $\bm\beta\ne\bm0$.

    The weighted sum of the selected inequalities is
    \[
        \bm\beta^{\top}\bm\theta\le\ell(\bm x),
        \qquad
        \ell(\bm x)=\sum_{a\in F}w_aq_a(\bm x)
                         +\sum_{j\in J}\mu_j\ell_j(\bm x).
    \]
    It is valid for $\mathcal P^+$ and active at the candidate.
    Equality in the fiber at $\bar{\bm x}$ forces every selected inequality to be active, and full rank makes $\bar{\bm\theta}$ its unique solution.
    Hence $\bar{\bm\theta}$ uniquely maximizes $\bm\beta^{\top}\bm\theta$ over $\mathcal P^+(\bar{\bm x})$.
    Because $\conv(\mathcal S)\subseteq\mathcal P^+$, the attained identity \eqref{eq:app-hull-fiber-value} gives exactly the same zero-violation equivalence as in the preceding proof.
    Finally, a vertex of $\mathcal P^+$ with fractional leader coordinates cannot belong to $\conv(\mathcal S)$, since $\mathcal S\subseteq\mathcal P^+$ and its leader coordinates are binary.
\end{proofof}

\subsection{Complexity of separation}\label{appx:proof-separation-hardness}

\begin{proofof}{\autoref{prop:separation-hardness}}
    Assume rational data with binary encoding.
    Let $s_1,\ldots,s_n$ be positive integer weights.
    The \textsc{Partition} problem~\cite{Karp1972} asks whether there exists a subset $I\subseteq\{1,\ldots,n\}$ such that
    \[
        \sum_{i\in I}s_i=\sum_{i\notin I}s_i.
    \]
    Write $W=\sum_{i=1}^n s_i$.
    Set $\mathcal X=\{0,1\}^n$ and give the two parallel arcs costs
    \[
        q_{a_1}(\bm x)=2\bm s^{\top}\bm x-W,
        \qquad q_{a_2}(\bm x)=0.
    \]
    The point $\bar{\bm x}=\tfrac12\bm1$, $\bar\theta_r=\bar\theta_t=0$ belongs to $\mathcal P$, with both arc inequalities active.
    Its root-supply function is $\phi(\bm x)=\min\{2\bm s^{\top}\bm x-W,0\}$, so $\hat\phi\le0$.
    The candidate belongs to $\setR_r$ if and only if $\hat\phi(\bar{\bm x})=0$.

    If $\bm x^*\in\{0,1\}^n$ satisfies $2\bm s^{\top}\bm x^*=W$, then its complement also does, and both have follower value zero.
    Their average is $\bar{\bm x}$, so \eqref{eq:app-envelope-lp} gives $\hat\phi(\bar{\bm x})=0$.
    Conversely, suppose this envelope value is zero and take attaining weights in \eqref{eq:app-envelope-lp}.
    Since all integer follower values are nonpositive, every positively weighted assignment satisfies $\phi(\bm x^k)=0$, hence $2\bm s^{\top}\bm x^k\ge W$.
    Their weighted average satisfies
    \[
        \sum_{k=1}^K\alpha_k\,2\bm s^{\top}\bm x^k
        =2\bm s^{\top}\bar{\bm x}=W.
    \]
    Every positively weighted assignment therefore satisfies equality and solves the \textsc{Partition} instance with weights $(s_1,\ldots,s_n)$.
    The construction has polynomial encoding length, and \textsc{Partition} is a well-known NP-complete problem, proving NP-hardness of membership under the stated restrictions.
    An exact root-cut separator decides this membership by reporting whether a violated cut exists, and exact evaluation of $\hat \phi$ decides it by comparison with zero.
    Both tasks are therefore NP-hard as well.
\end{proofof}

\end{document}